\documentclass[12pt,reqno]{article}
\usepackage
[colorlinks=true, linkcolor=webgreen, filecolor=webbrown,
citecolor=webgreen]{hyperref}
\def\N{{\mathbb N}}

\usepackage{amsmath}
\usepackage{amssymb}
\usepackage{color}
\usepackage{epsf}

\begin{document}
\makeatletter

\begin{center}
\epsfxsize=10in
\end{center}

\def\endofproofmark{$\Box$}

\begin{center}
\vskip 1cm {\LARGE\bf Recurrence Divisibility Tests} \vskip 1cm
\vspace{10mm}

{\large Mehdi Hassani}

\vskip .5cm
Department of Mathematics\\
Institute for Advanced
Studies in Basic Sciences\\
Zanjan, Iran\\
\href{mailto:mhassani@iasbs.ac.ir}{\tt mhassani@iasbs.ac.ir}
\end{center}

\newtheorem{thm}{Theorem}
\newtheorem{prop}{Proposition}
\newtheorem{lemma}{Lemma}
\newtheorem{cor}{Corollary}
\newtheorem{prob}{Note and Problem}
\def\frameqed{\framebox(5.2,6.2){}}
\def\deshqed{\dashbox{2.71}(3.5,9.00){}}
\def\ruleqed{\rule{5.25\unitlength}{9.75\unitlength}}
\def\myqed{\rule{8.00\unitlength}{12.00\unitlength}}
\def\qed{\hbox{\hskip 6pt\vrule width 7pt height11pt depth1pt\hskip 3pt}
\bigskip}
\newenvironment{proof}{\trivlist\item[\hskip\labelsep{\bf Proof}:]}{\hfill
 $\frameqed$ \endtrivlist}
\newcommand{\COM}[2]{{#1\choose#2}}

\thispagestyle{empty} \null \addtolength{\textheight}{1cm}

\begin{abstract}
In this note, we are going to introduce some recurrence
divisibility tests for all primes except than 2 and 5.
\end{abstract}

\bigskip
\hrule
\bigskip


\bigskip
\hrule
\bigskip
\vskip .2cm\hspace{-8mm} For $n\in\N$, we define $ud(n)=$ unite
digit of $n$. As we know, we have the following divisibility tests
for 2 and 5:
$$
2|n\Leftrightarrow 2|ud(n),
$$
and
$$
5|n\Leftrightarrow ud(n)\in\{0,5\}.
$$
In this note, we are going to find some \textit{recurrence
divisibility tests} for other primes except than 2 and 5. As soon,
you will understand our mind by ``recurrence''.
\begin{thm} For every prime $p\neq 2, 5$, there exists unique
$t(p)\in\mathbb{Z}_p$, in which
$$
p|n\Leftrightarrow
p\Big|\Big\lfloor\frac{n}{10}\Big\rfloor-t(p)ud(n).
$$
Also, we have
$$
t(p)=\left\{ \begin{array}{ll}
\frac{p-1}{10} & {\rm if~} ud(p)=1,\\
\frac{7p-1}{10} & {\rm if~} ud(p)=3,\\
\frac{3p-1}{10} & {\rm if~} ud(p)=7,\\
\frac{9p-1}{10} & {\rm if~} ud(p)=9.
\end{array}
\right.
$$
\end{thm}
\begin{proof} Suppose $p\neq 2, 5$ is a prime. Since
gcd$(10,p)=1$, the diophantine equation $10t+1=kp$ with $0\leq
t<p$, has an unique solution module $p$ and clearly $kt>0$. Let
$t=t(p)$. So, $p|10t(p)+1$ and $1\leq t(p)\leq p-1$.\\
Now, suppose $p|n$, or consequently $p|t(p)n$. Also, we have
$p|10t(p)+1$. Therefore, by considering
$ud(n)=n-10\lfloor\frac{n}{10}\rfloor$, we have
$p|\lfloor\frac{n}{10}\rfloor-t(p)ud(n)$.\\
Also, suppose $p|\lfloor\frac{n}{10}\rfloor-t(p)ud(n)$ or
$p|(10t(p)+1)\lfloor\frac{n}{10}\rfloor-t(p)n$. Since
$p|10t(p)+1$, we obtain $p|t(p)n$ and since $1\leq t(p)\leq p-1$,
we have gcd$(p,t(p))=1$, therefore $p|n$.\\
Now, we compute the value of $t(p)$. To do this, we note that
since $p\neq 2, 5$, we have $ud(p)\in\{1,3,7,9\}$. If we let
$k=k(p)$, we have the relation $10t(p)+1=k(p)p$ with $1\leq
t(p)\leq p-1$. So, $ud(pk(p))=1$ always holds. Also, the condition
$1\leq t(p)\leq p-1$ yield that $\frac{9}{p}\leq k(p)\leq
10-\frac{11}{p}$. So, for all primes $p\neq 2, 5$ we have $1\leq
k(p)\leq 9$. If $ud(p)=1$, then since $ud(pk(p))=1$ and $1\leq
k(p)\leq 9$, we obtain $k(p)=1$ and so $t(p)=\frac{p-1}{10}$.
Other cases have similar reasons. This completes the proof.
\end{proof}
\textbf{Note 1.} We called this divisibility test, ``recurrence'',
because for all primes $p\neq 2, 5$ we have
$$
\lfloor\log_{10}n\rfloor=
1+\Big\lfloor\log_{10}\Big(\Big\lfloor\frac{n}{10}\Big\rfloor-t(p)ud(n)\Big)\Big\rfloor.
$$
It means, this test works by cancelling digits of given number for
check, 1 digit by 1 digit.\\\\
\textbf{Note 2.} Same divisibility tests can be obtain which work
by cancelling $m$ digit by $m$ digit.\\\\
\textbf{Note 3.} We can yield same tests for all $n\in\N$ with $2\nmid n$ and $5\nmid n$.\\\\
\textbf{Acknowledgements.} I deem my duty to thank Mr. H. Osanloo,
who mentioned me recurrence divisibility tests in some special
cases.

\end{document}